\documentclass[12pt]{article}%
\usepackage{amsmath}
\usepackage{amsfonts}
\usepackage{amssymb}
\usepackage{graphicx}
\usepackage{color}%
\providecommand{\U}[1]{\protect\rule{.1in}{.1in}}
\makeindex

\begin{document}

\title{\textbf{Projection Methods: An Annotated Bibliography of Books and Reviews }}
\author{Yair Censor$^{1}$ and Andrzej Cegielski$^{2}$\medskip\ \\$^{1}${Department of Mathematics, University of Haifa,}\\Mt. Carmel, Haifa 3498838, Israel\medskip\\$^{2}${Faculty of Mathematics, Computer Science and}\\Econometrics, University of Zielona G\'{o}ra, \\ul. Szafrana 4a, 65-514 Zielona G\'{o}ra, Poland}
\date{July 30, 2014. Revised: September 5, 2014.}
\maketitle
\tableofcontents

\begin{abstract}
Projections onto sets are used in a wide variety of methods in optimization
theory but not every method that uses projections really belongs to the class
of projection methods as we mean it here. Here \textit{projection methods} are
iterative algorithms that use projections onto sets while relying on the
general principle that when a family of (usually closed and convex) sets is
present then projections (or approximate projections) onto the given
individual sets are easier to perform than projections onto other sets
(intersections, image sets under some transformation, etc.) that are derived
from the given family of individual sets. Projection methods\textit{ }employ
projections (or approximate projections) onto convex sets in various ways.
They may use different kinds of projections and, sometimes, even use different
projections within the same algorithm. They serve to solve a variety of
problems which are either of the feasibility or the optimization types. They
have different algorithmic structures, of which some are particularly suitable
for parallel computing, and they demonstrate nice convergence properties
and/or good initial behavioral patterns. This class of algorithms has
witnessed great progress in recent years and its member algorithms have been
applied with success to many scientific, technological, and mathematical
problems. This annotated bibliography includes books and review papers on, or
related to, projection methods that we know about, use, and like. If you know
of books or review papers that should be added to this list please contact us.

\end{abstract}

\bigskip

\textbf{2010 Mathematics Subject Classification (MSC)}: 65K10, 90C25

\section{Introduction}

\textbf{Projection methods}. Projections onto sets are used in a wide variety
of methods in optimization theory but not every method that uses projections
really belongs to the class of projection methods as we mean it here. Here
\textit{projection methods} are iterative algorithms that use projections onto
sets while relying on the general principle that when a family of (usually
closed and convex) sets is present then projections (or approximate
projections) onto the given individual sets are easier to perform than
projections onto other sets (intersections, image sets under some
transformation, etc.) that are derived from the given family of individual sets.

A projection algorithm reaches its goal, related to the whole family of sets,
by performing projections onto the individual sets. Projection
algorithms\textit{\ }employ projections (or approximate projections) onto
convex sets in various ways. They may use different kinds of projections and,
sometimes, even use different projections within the same algorithm. They
serve to solve a variety of problems which are either of the feasibility or
the optimization types. They have different algorithmic structures, of which
some are particularly suitable for parallel computing, and they demonstrate
nice convergence properties and/or good initial behavioral patterns in some
significant fields of applications.

Apart from theoretical interest, the main advantage of projection methods,
which makes them successful in real-world applications, is computational. They
commonly have the ability to handle huge-size problems of dimensions beyond
which other, more sophisticated currently available, methods start to stutter
or cease to be efficient. This is so because the building bricks of a
projection algorithm are the projections onto the given individual sets
(assumed and actually easy to perform) and the algorithmic structures are
either sequential or simultaneous or in-between, such as in the
block-iterative projection (BIP) methods or in the more recent
string-averaging projection (SAP) methods. An advantage of projection methods
is that they work with initial data and do not require transformation of, or
other operations on, huge and sparse matrices describing the problem.

\textbf{Purpose of the paper}. We present here an annotated bibliography of
books and review papers on projection methods. This should be helpful to
researchers, veterans or newcomers, in the field by directing them to the many
existing resources. The vast amount of research papers in the field of
projection methods\ makes it sometimes difficult to master even within a
specific sub-area. On the other hand, projection methods send branches both
into fields of applications wherein real-world problems are solved and into
theoretical areas in mathematics such as, but not only, fixed point theory and
variational inequalities. Researchers in each of these, seemingly
perpendicular, directions might wish to consult the bibliography presented
here and check the annotated items' references lists for further information.

We emphasize that we had no intention to create a bibliography of the research
itself in the field of projection methods. This would have been a humongous
job beyond our plans. So, we stick with books, book chapters and review
papers, except for a few glitches that made us unable to refrain from
mentioning a (very) few articles which are not strictly review papers but
contain beside their research contents also some worthwhile review-like material.

\textbf{An apology}. Oversight and lack of knowledge are human traits which we
are not innocent of. Therefore, we apologize for omissions and other
negligence and lacunas in this paper. We kindly ask our readers to communicate
to us any additional items and informations that fit the structure and spirit
of the paper and we will gladly consider those for inclusion in future
revisions, extensions and updates of the paper. Such updates will be posted on www.arxiv.org.

\textbf{Organization of the paper}. After paying tribute to some of the early
workers who laid the foundations (Section \ref{sect:early}) we cite and
annotate books in Section \ref{sect:books}, followed by annotated citations of
review papers in Section \ref{sect:papers}. We go alphabetically through the
cited items in each section.

\section{Early beginnings: von Neumann, Kaczmarz and Cimmino\label{sect:early}%
}

In telegraphic language we recognize von Neumann, Kaczmarz and Cimmino as the
forefathers. John von Neumann's 1933 \cite{neumann} method of alternating
projections (MAP) is a projection method for finding the projection of a given
point onto the intersection of two subspaces in Hilbert space. Stefan
Kaczmarz, in a three pages paper published in 1937 \cite{kacmarcz},
(posthumous translation into English in \cite{kacmarcz-eng93}) presented a
sequential projections method for solving a (consistent) system of linear
equations. Historical information about his life can be found in the papers of
Maligranda \cite{maligranda} and Parks \cite{parks93}.

Gianfranco Cimmino proposed in \cite{cimmino}, published in 1938, a
simultaneous projection method for the same problem in which one projects the
current iterate simultaneously on all hyperplanes, representing the linear
equations of the system, and then takes a convex combination to form the next
iterate. A historical account of Cimmino's work was published by Benzi
\cite{benzi05}.

In 1954 Agmon \cite{agmon} and Motzkin and Schoenberg \cite{mot-schoen54}
generalized the sequential projections method from hyperplanes to half-spaces,
and then Eremin \cite{eremin65} in 1965, Bregman \cite{Bregman} in 1965, and
Gubin, Polyak and Raik in 1967 \cite{GPR67} generalized it farther to convex
sets. In 1970 Auslender \cite[page 78]{auslender} generalized Cimmino's
simultaneous projection method to convex sets. Arriving from a different
perspective, one must also mention here the seminal papers of Amemiya and Ando
\cite{amemiya-ando} and of Halperin \cite{halperin} that discuss products of
contractions and products of projection operators, respectively, in Hilbert space.

A significant early contribution was the work of Lev Bregman who introduced a
new \textquotedblleft distance\textquotedblright\ that gave rise to new
\textquotedblleft projections\textquotedblright\ that include as special cases
the Euclidean distance and projection as well as the Kullback-Leibler
\textquotedblleft distance\textquotedblright\ and entropic projections.
Bregman's 1967 seminal paper \cite{bregman67}, based on his doctoral
dissertation, had absolutely no follow-up in the literature until the
appearance, 14 years later, of \cite{cl81}. Nowadays, Bregman projections and
projection algorithms that employ them are abundant in the literature.

These were the early beginnings that paved the way for the subsequent
\textquotedblleft explosion\textquotedblright\ of research in this field that
continues to this day and covers many aspects. These include, but are not
limited to, developments of new algorithmic structures for projection methods,
usage of different types of projections, application of projection methods to
new types of feasibility, optimization, or variational inequalities problems,
investigations of the above in various spaces, branching into fixed point
theory and other mathematical areas, and using projection methods in
significant real-world problems with real data sets of humongous dimensions,
and more.

The interface between projection methods and significant real-world problems
was and still is a fertilizer for both. In 1970 Gordon, Bender and Herman
\cite{gbh71} published an algebraic reconstruction technique (ART) for
three-dimensional electron microscopy and x-ray photography. It was recognized
later on that their ART is an independent rediscovery of Kaczmarz's sequential
projections method for solving a system of linear equations. The first CT
(computerized tomography) scanner by EMI (Electric \& Musical Industries Ltd.,
London, England, UK), invented by G.N. Hounsfield \cite{hounsfield}, used a
variant of ART, i.e., of Kaczmarz's sequential projections method. For his
pioneering invention, Hounsfield shared a Nobel prize with A.M. Cormack in
1979. Read more on this field in \cite{GTH}.

\section{Books on projection methods\label{sect:books}}

\cite{auslender} \textbf{A. Auslender, \textit{Optimisation: M\'{e}thodes
Num\'{e}riques}, Masson, Paris, France, 1976 (in French).\bigskip}

This book on optimization is one of the earliest books to include a subsection
on the convex feasibility problem in the finite-dimensional Euclidean space
$R^{n}$ (in Chapter V) and to discuss orthogonal (least Euclidean distance)
projections onto convex sets. To our knowledge, the algorithm in Equations
(1.6)--(1.7) on page 78 therein is the first published version of Cimmino's
simultaneous projections algorithm extended to handle convex
sets.\textbf{\bigskip}

\cite{bauschke-conf} \textbf{H.H. Bauschke, R.S. Burachik, P.L. Combettes, V.
Elser, D.R. Luke and H. Wolkowicz (Editors), \textit{Fixed-Point Algorithms
for Inverse Problems in Science and Engineering}, Springer, 2011.\bigskip}

This book presents recent work in variational and numerical analysis. The
contributions provide state-of-the-art theory and practice in first-order
fixed-point algorithms, identify emerging problems driven by applications, and
discuss new approaches for their solution. The book is a compendium of topics
explored at the Banff International Research Station \textquotedblleft
Interdisciplinary Workshop on Fixed-Point Algorithms for Inverse Problems in
Science and Engineering\textquotedblright\ in November 2009. The link to
projection methods is due to the mathematical topics discussed which include:
Bregman distances, feasibility problems, the common fixed point problem, the
Douglas--Rachford algorithm, monotone operators, proximal splitting methods,
nonexpansive mappings and more, all of which are highly relevant to projection
methods.\textbf{\bigskip}

\cite{BC11} \textbf{H.H. Bauschke and P.L. Combettes, \textit{Convex Analysis
and Monotone Operator Theory in Hilbert Spaces}, Springer, New York, NY, USA,
2011.\bigskip}

This outstanding book presents relations among convex analysis, monotone
operators and nonexpansive operators. The fact that the metric projection
operator is both monotone and nonexpansive in not only reason for including
this book here. More important is the fact that the metric projection, as a
firmly nonexpansive operator, is the resolvent of a maximal monotone operator.
Moreover, finding a zero of a monotone operator is equivalent to finding a
fixed point of its resolvent. The authors present a wide spectrum of
properties of all these operators and apply these properties to many methods
which use projections onto closed convex subsets.\textbf{\bigskip}

\cite{bg03} \textbf{A. Ben-Israel and T. Greville, \textit{Generalized
Inverses: Theory and Applications}, 2nd Edition, Springer-Verlag, New York,
NY, USA, 2003.\bigskip}

This fundamental work on generalized inverses discusses orthogonal projections
and projectors (Chapter 2, Section 7) and projectors associated with
essentially strictly convex norms (Chapter 3, Section 7) and relates them to
generalized inverses.\textbf{\bigskip}

\cite{berinde} \textbf{V. Berinde, \textit{Iterative Approximation of Fixed
Points}, Lecture Notes in Mathematics, Vol. 1912, Springer-Verlag, Berlin,
Heidelberg, Germany, 2007.\bigskip}

This book deals with convergence and with stability of fixed point iterative
procedures in Banach spaces. Most of these procedures are defined by Pickard
iteration, Krasnosel'ski\u{\i} iteration, Mann iteration, Ishikawa iteration,
viscosity approximation, or their modifications. The convergence theorems
presented in the book can be in particular applied to various projection
methods.\textbf{\bigskip}

\cite{BO75} \textbf{E. Blum and W. Oettli, \textit{Mathematische Optimierung:
Grundlagen und Verfahren}, Springer-Verlag, Berlin, Heidelberg, Germany,
1975.\bigskip}

This 1975 German language book on mathematical optimization contains a chapter
(Kapitel 6) on projections and contractive methods. In it, the projected
gradient minimization method is studied (under the name \textquotedblleft the
method of Uzawa\textquotedblright) and Fej\'{e}r-contractive mappings are
handled. Those are then used to describe an iterative method of Eremin and
Mazurov.\bigskip

\cite{borwein-zhu-book} \textbf{J.M. Borwein and Q.J. Zhu, \textit{Techniques
of Variational Analysis}, Springer-Verlag, New York, NY, USA, 2005. Paperback
2010.}\bigskip

Variational arguments are classical techniques whose use can be traced back to
the early development of the calculus of variations and further. The book
discusses various forms of variational principles in Chapter 2 and then
discusses applications of variational techniques in different areas in
Chapters 3--7. Section 4.5 in Chapter 4 is devoted to Convex Feasibility
Problems since they can be viewed as, and handled by, techniques of
variational analysis.\bigskip

\cite{brezinski-book} \textbf{C. Brezinski, \textit{Projection Methods for
Systems of Equations}, Elsevier Science Publishers, Amsterdam, The
Netherlands, 1997.\bigskip}

The main part of this book is devoted to iterative solution of systems of
linear equations in a Euclidean space. Most of the methods described in the
book employ projections (orthogonal or oblique) onto a sequence of
hyperplanes. The author gives a general model and shows that many projection
methods known from the literature are special cases of it. The book also
describes how to accelerate the convergence of these methods.\textbf{\bigskip}

\cite{pink-book} \textbf{D. Butnariu, Y. Censor and S. Reich (Editors),
\textit{Inherently Parallel Algorithms in Feasibility and Optimization and
Their Applications}, Elsevier Science Publishers, Amsterdam, The Netherlands,
2001.\bigskip}

This 504 pages book contains papers presented at the \textquotedblleft
Research Workshop on Inherently Parallel Algorithms in Feasibility and
Optimization and Their Applications\textquotedblright, held March 13--16,
2000, jointly at the University of Haifa and the Technion in Haifa, Israel.
Most of the 27 papers in it are about projection methods or closely related to
them and the focus is on parallel feasibility and optimization algorithms and
their applications, in particular on inherently parallel algorithms. By this
term one means algorithms which are logically (i.e., in their mathematical
formulations) parallel, not just parallelizable under some conditions, such as
when the underlying problem is decomposable in a certain
manner.\textbf{\bigskip\ }

\cite{byrnebook} \textbf{C.L. Byrne, \textit{Applied Iterative Methods}, AK
Peters, Wellsely, MA, USA, 2008.\bigskip}

The author describes in this book a huge number of iterative methods for
several optimization problems in a Euclidean space. The main part of the book
is devoted to projection methods. The author presents many methods and their
variants known from the literature, as well as their
applications.\textbf{\bigskip}

\cite{byrne14book} \textbf{C.L. Byrne, \textit{Iterative Optimization in
Inverse Problems}, Chapman and Hall/CRC Press, Boca Raton, FL, USA,
2014.\bigskip}

The books is devoted to iterative methods for inverse problems in a Euclidean
space, in which projection methods play an important role. In Chapters 6 and 7
the author presents properties of several classes of operators containing the
metric projections onto closed convex subsets. In Chapters 8 and 9 these
properties are applied to projection methods for the convex feasibility
problem and for the split feasibility problem.\textbf{\bigskip}

\cite{CEG12} \textbf{A. Cegielski, \textit{Iterative Methods for Fixed Point
Problems in Hilbert Spaces}, Lecture Notes in Mathematics, Vol. 2057,
Springer-Verlag, Berlin, Heidelberg, Germany, 2012.\bigskip}

The book deals with iteration methods for solving fixed points problems in a
Hilbert space. The presentation consolidates many methods which apply
nonexpansive and quasi-nonexpansive operators. The author gives a wide
spectrum of properties of several classes of operators: quasi-nonexpansive,
strongly quasi-nonexpansive, nonexpansive and averaged. All these classes
contain the metric projections onto closed convex subsets and are closed under
convex combination and composition. The properties of these classes, as well
as some general convergence theorems, enable to prove the convergence for
members of a large class of projection methods for solving fixed point
problems.\textbf{\bigskip}

\cite{CZ97} \textbf{Y. Censor and S.A. Zenios, \textit{Parallel Optimization:
Theory, Algorithms, and Applications}, Oxford University Press, New York, NY,
USA, 1997.\bigskip}

This book presents Bregman distances and Bregman projections along with
Csiszar $\varphi$-divergences (Chapter 2). Many projection methods are then
described along with their convergence analyses. Chapter 5 is about iterative
methods for convex feasibility problems, all of which are projection methods.
Chapter 6 deals with iterative algorithms for linearly constrained
optimization problems, all of which are again from the family of projection
methods and embedded in the more general framework of Bregman projections (of
which orthogonal projections are a special case). The third part of the book
is devoted to applications and has chapters for matrix estimation problems,
image reconstruction from projections, the inverse problem of radiation
therapy treatment planning, multicommodity network flow problems and planning
under uncertainty.\textbf{\bigskip}

\cite{chinneck-book} \textbf{J.W. Chinneck, \textit{Feasibility and
Infeasibility in Optimization: Algorithms and Computational Methods},
Springer, New York, NY, USA, 2007.\bigskip}

Part I of this book is entitled: \textquotedblleft Seeking
Feasibility\textquotedblright, Part II is entitled: \textquotedblleft
Analyzing Infeasibility\textquotedblright\textbf{ }and Part III is devoted to
applications. With emphasis on the computational aspects and without shying
from heuristic algorithmic modifications that extend the algorithms'
efficiency, the book summarizes the state-of-the-art at the interface of
optimization and feasibility. Anyone interested in feasibility-seeking methods
that work or anyone looking for algorithms that have not yet been fully
explored mathematically will find here interesting materials.\textbf{\bigskip
\ }

\cite{deutsch01} \textbf{F. Deutsch, \textit{Best Approximation in Inner
Product Spaces}, Springer-Verlag, New York, NY, USA, 2001.\bigskip}

Frank Deutsch's book is a fundamental resource for all who study the theory of
projection methods in Hilbert spaces. The core of the book is the metric
projection onto a closed convex subset and its properties. Chapters 9-10 are
devoted to projection methods. The author applies the properties of the metric
projection to von Neumann's alternating projection method and to Halperin's
cyclic projection method for solving the best approximation problem for a
finite family of closed subspaces. Further, the Dykstra algorithm is presented
for solving the best approximation problem for a finite family of closed
convex subsets. The book contains many examples, where the presented methods
are applied to linear equations, linear inequalities, isotone regression,
convex regression and to the shape-preserving interpolation. Every chapter in
the book ends with interesting exercises.\textbf{\bigskip}

\cite{eremin}\textbf{ I.I. Eremin, \textit{Theory of Linear Optimization}, VSP
-- International Science Publishers, Zeist, The Netherlands, 2002.\bigskip}

This monograph is dedicated to the basic component of the theory of linear
optimization: systems of linear inequalities. Chapter 6 is devoted to methods
of projection in linear programming and deals with Fej\'{e}r mappings and
processes in that context.\textbf{ }This is a topic on which the author has
published research papers, as mentioned in Section \ref{sect:early}.\bigskip

\cite{ER11} \textbf{R. Escalante and M. Raydan, \textit{Alternating Projection
Methods}, Society for Industrial and Applied Mathematics (SIAM), Philadelphia,
PA, USA, 2011.\bigskip}

Starting from von Neumann's 1933 method of alternating projections (MAP) for
finding the projection of a given point onto the intersection of two subspaces
in Hilbert space, Escalante and Raydan meticulously make their way to
row-action methods, Bregman projections, the convex feasibility problem,
Dykstra's algorithm for the best approximation problem (BAP) and matrix
problems. Needless to say that projection methods are the backbone of this
book which appeared in SIAM's series on \textquotedblleft Fundamentals of
Algorithms\textquotedblright.\textbf{\bigskip}

\cite{PF03} \textbf{F. Facchinei and J.-S. Pang, \textit{Finite-Dimensional
Variational Inequalities and Complementarity Problems}, Volume I and Volume
II, Springer-Verlag, New York, NY, USA, 2003.\bigskip}

This over 1200 pages two volumes book is a comprehensive compendium on theory
and methods for solving variational inequalities in Euclidean spaces. In
particular, Chapter 12 is devoted to iterative methods for monotone problems
where the metric projections play an important role. The authors present many
methods: fixed point iterations, extragradient methods, hyperplane projection
methods, regularization algorithms, all of which employ metric projections
(orthogonal or oblique) in each iteration. The rest of this chapter is devoted
to more general methods: proximal point methods and splitting
methods.\textbf{\bigskip}

\cite{galantai} \textbf{A. Gal\'{a}ntai, \textit{Projectors and Projection
Methods}, Kluwer Academic Publishers, Dordrecht, The Netherlands,
2004.\bigskip}

With idempotent matrices as projectors, this book is mostly linear algebra
oriented in its presentation and the techniques that it uses. As such it helps
to look at projection methods from that perspective which is complementary to
the analysis perspective exercised by other books mentioned in this paper.
Chapter 4 on iterative projection methods for linear algebraic systems
includes descriptions of the methods of Kaczmarz, Cimmino, Altman, Gastinel,
Householder and Bauer. The 369 items long references list is a treasure chest
for anyone interested in the field.\textbf{\bigskip}

\cite{gastinel} \textbf{N. Gastinel, \textit{Linear Numerical Analysis},
Hermann, Paris, France, 1970. Translated from the original French text:
\textit{Analyse Num\'{e}rique Lin\'{e}aire}, Hermann, Paris, France,
1966.\bigskip}

Limited to linear numerical analysis, Gastinel's 1970 book contains a chapter
(Chapter 5) dedicated to \textquotedblleft indirect methods\textquotedblright%
\ for solving linear systems. Orthogonal (least Euclidean distance)
projections onto hyperplanes are considered and Kaczmarz's method and
Cimmino's method for linear systems are studied. We think that this is the
first appearance of these methods in a textbook, although Kaczmarz's method
and its applications were described by Tompkins in Chapter 18 of an edited
book \cite{beckenbach} published in 1956.\textbf{\bigskip}

\cite{GTH} \textbf{G.T. Herman, \textit{Fundamentals of Computerized
Tomography: Image Reconstruction from Projections}, Springer-Verlag, London,
UK, 2nd Edition, 2009.\bigskip}

This is a revised and updated version of the successful 1980 first edition.
Warning: the term \textquotedblleft projections\textquotedblright\ in
\textquotedblleft image reconstruction from projections\textquotedblright\ has
a different meaning than \textquotedblleft projections onto
sets\textquotedblright, as used elsewhere in this paper. The book is devoted
to the fundamentals of the field of image reconstruction from
projections\textbf{ }with particular emphasis on the computational and
mathematical procedures underlying the data collection, image reconstruction,
and image display in the practice of computerized tomography. It is written
from the point of view of the practitioner: points of implementation and
application are carefully discussed and illustrated. The major emphasis of the
book is on reconstruction methods; these are thoroughly surveyed. Chapter 11
on algebraic reconstruction techniques and Chapter 12 on quadratic
optimization methods are related to projection methods as they are used in
this field.\textbf{\bigskip}

\cite{householder64} \textbf{A.S. Householder, \textit{The Theory of Matrices
in Numerical Analysis}, Dover Publications, Inc., New York, NY, USA,
1975.\bigskip}

\textbf{ }Originally published in 1964,\textbf{ }like good wine, aging does
not diminish the value and beauty of this book. The Preface opens with:
\textquotedblleft This book represents an effort to select and present certain
aspects of the theory of matrices that are most useful in developing and
appraising computational methods for solving systems of linear
equations...\textquotedblright. We think that the book lives up to this,
particularly since Chapter 4 on \textquotedblleft The Solution of Linear
Systems: Methods of Successive Approximation\textquotedblright\ includes the
Subsection 4.2 on \textquotedblleft Methods of Projection\textquotedblright%
\ which makes this book eligible for inclusion here. Besides Kaczmarz's method
it mentions also an 1951 projection method due to de la Garza.\bigskip

\cite{konnov}\textbf{ I. Konnov, \textit{Combined Relaxation Methods for
Variational Inequalities}, Lecture Notes in Economics and Mathematical
Systems, Vol. 495, Springer-Verlag, Berlin, Heidelberg, Germany,
2001.}\bigskip

The combined relaxation (CR) methods of Konnov do not fall into our definition
of \textquotedblleft projection methods\textquotedblright. On the contrary,
they use a projection step onto a set within a more general two-level
algorithmic structure. However, the activity in CR methods to replace
projections onto a set by projections onto separating hyperplanes has
something in common with similar activities in \textquotedblleft projection
methods\textquotedblright\ as we mean them here. Besides, the careful
treatment of such algorithms for variational inequalities and the research
possibilities opened by looking at those from the point of view of projection
methods tilted our decision towards including the book here.\bigskip

\cite{kurpel} \textbf{N.S. Kurpel', \textit{Projection-Iterative Methods for
Solution of Operator Equations}, Translations of Mathematical Monographs, Vol.
46, American Mathematical Society, Providence, RI, USA, 1976.\bigskip}

This is a 196 pages monograph with 238 bibliographic items many of which are
not cited anymore in current days literature. Iterative methods for solution
of operator equation are studied here in a general setting. Abstract metric
and normed spaces (in Chapter 1) as well as Banach spaces (in Chapter 2) are
the background against which iterative methods, employing general algorithmic
operators, are investigated. Specialization of those general algorithmic
operators to projection operators leads to projection methods.
\textbf{\bigskip}

\cite{p12} \textbf{C. Popa, \textit{Projection Algorithms - Classical Results
and Developments: Applications to Image Reconstruction}, Lambert Academic
Publishing - AV Akademikerverlag GmbH \& Co. KG, Saarbr\"{u}cken, Germany,
2012.\bigskip}

As a textbook, this book provides a short and useful introduction into the
field of projection based solvers. It is at the same time also a research
monograph that describes some of the author's research results in the field.
Set in the finite-dimensional Euclidean space context, the theoretical
discussions provide ideas for further developments and research.\textbf{
}Besides\textbf{ }introductory material, the book includes chapters on
extensions to inconsistent least squares problems, oblique and generalized
oblique projections, constraining strategies, and some special
projection-based algorithms.\textbf{\bigskip}

\cite{saad03} \textbf{Y. Saad, \textit{Iterative Methods for Sparse Linear
Systems}, Second Edition, Society for Industrial and Applied Mathematics
(SIAM), Philadelphia, PA, USA, 2003.\bigskip}

The book, freely available from the author's homepage on the Internet,
declares its intention to provide up-to-date coverage of iterative methods for
solving large sparse linear systems. It focuses on practical methods that work
for general sparse matrices rather than for any specific class of problems.
Although the very definition of \textquotedblleft projection
methods\textquotedblright, in Chapter 5 (\textquotedblleft Projection
Methods\textquotedblright) is not identical with what we mean by this term
here, the relations and connections between them cannot be mistaken.
\textbf{\bigskip}

\cite{itor} \textbf{H.D. Scolnik, A.R. De Pierro, N.E. Echebest and M.T.
Guardarucci (Guest Editors), \textit{Special Issue of International
Transactions in Operational Research}, Volume 16, Issue 4, July 2009.\bigskip}

This special issue is devoted to projection methods. From block-iterative
algorithms, through the string-averaging method for sparse common fixed-point
problems, and from convergence of the method of alternating projections
through perturbation-resilient block-iterative projection methods with
application to image reconstruction from projections -- this special issue is
a focused source for novel ideas and for literature coverage of the subject.
It is available at:
http://onlinelibrary.wiley.com/doi/10.1111/itor.2009.16.issue-4/issuetoc.\textbf{\bigskip
}

\cite{stark} \textbf{H. Stark and Y. Yang, \textit{Vector Space Projections: A
Numerical Approach to Signal and Image Processing, Neural Nets, and Optics},
John Wiley \& Sons, Inc. New York, NY, USA, 1998.\bigskip}

Using the term \textquotedblleft vector-space projections\textquotedblright,
this book presents a nice blend of theory (Chapters 1--5) and applications
(Chapters 6--9) of the \textquotedblleft projections onto convex
sets\textquotedblright\ (POCS) methods which is the name adopted mostly by the
engineering community for projection methods. The book reflects the growing
interest in the application of these methods to problem solving in science and
engineering. It brings together material previously scattered in disparate
papers, book chapters, and articles, and offers a systematic treatment of
vector space projections. \textbf{\bigskip}

\section{Review papers on projection methods\label{sect:papers}}

\cite{bb96} \textbf{H.H. Bauschke and J.M. Borwein, On projection algorithms
for solving convex feasibility problems,\ \textit{SIAM Review} 38 (1996),
367--426.\bigskip}

This review paper, based on the first author's Ph.D. work, has not lost its
vitality to this day. With 109 items in its bibliography and a subject index,
it is a treasure of knowledge central to the field. The abstract says:
\textquotedblleft Due to their extraordinary utility and broad applicability
in many areas of classical mathematics and modern physical sciences (most
notably, computerized tomography), algorithms for solving convex feasibility
problems continue to receive great attention. To unify, generalize, and review
some of these algorithms, a very broad and flexible framework is investigated.
Several crucial new concepts which allow a systematic discussion of questions
on behavior in general Hilbert spaces and on the quality of convergence are
brought out. Numerous examples are given.\textquotedblright\textbf{\bigskip}

\cite{bbl97} \textbf{H.H. Bauschke, J.M. Borwein and A.S. Lewis, The method of
cyclic projections for closed convex sets in Hilbert space,
\textit{Contemporary Mathematics} 204 (1997), 1--38.\bigskip}

From the Abstract: \textquotedblleft Although in many applications [...] the
convex constraint sets do not necessarily intersect, the method of cyclic
projections is still employed. Results on the behaviour of the algorithm for
this general case are improved, unified, and reviewed. The analysis relies on
key concepts from convex analysis and the theory of nonexpansive
mappings.\textquotedblright\ Meticulously written, this paper contains a
wealth of material that both reviews and extends many important
results.\textbf{\bigskip\ }

\cite{swiss knives} \textbf{H.H. Bauschke and V.R. Koch, Projection methods:
Swiss army knives for solving feasibility and best approximation problems with
halfspaces, in: S. Reich and A. Zaslavski (Editors), Proceedings of the
workshop on Infinite Products of Operators and Their Applications, Haifa,
Israel, 2012, \textit{Contemprary Mathematics, }accepted for
publication.\bigskip}

Although geared toward solving the specific problem of automated design of
road alignments, this paper is written with an eye on reviewing the field. It
provides a selection of state-of-the-art projection methods, superiorization
algorithms, and best approximation algorithms. Various observations on the
algorithms and their relationships are given along with broad numerical
experiments introducing performance profiles for projection
methods.\textbf{\bigskip}

\cite{benzi05} \textbf{M. Benzi, Gianfranco Cimmino's contributions to
numerical mathematics, in: \textit{Atti del Seminario di Analisi Matematica},
Dipartimento di Matematica dell'Universita` di Bologna. Volume Speciale: Ciclo
di Conferenze in Memoria di Gianfranco Cimmino, Marzo-Aprile 2004, Tecnoprint,
Bologna, Italy (2005), pp. 87--109.\bigskip}

The abstract of this paper precisely describes its contents. It says:
\textquotedblleft Gianfranco Cimmino (1908--1989) authored several papers in
the field of numerical analysis, and particularly in the area of matrix
computations. His most important contribution in this field is the iterative
method for solving linear algebraic systems that bears his name, published in
1938. This paper reviews Cimmino's main contributions to numerical
mathematics, together with subsequent developments inspired by his work. Some
background information on Italian mathematics and on Mauro Picone's Istituto
Nazionale per le Applicazioni del Calcolo, where Cimmino's early numerical
work took place, is provided. The lasting importance of Cimmino's work in
various application areas is demonstrated by an analysis of citation patterns
in the broad technical and scientific literature.\textquotedblright%
\textbf{\bigskip}

\cite{bruck} \textbf{R.E. Bruck, On the random product of orthogonal
projections in Hilbert space II, \textit{Contemporary}}
\textbf{\textit{Mathematics} 513 (2010), 65--98.\bigskip}

This paper is mainly concerned with certain abstract properties of products of
linear orthogonal projections onto closed subspaces of a Hilbert space. We
mention it here because its first section presents a brief, but very readable
and informative, history of the study of convergence of infinite products of
such projections, as well as of their nonlinear counterparts.\textbf{\bigskip}

\cite{yc-siam-81} \textbf{Y. Censor, Row-action methods for huge and sparse
systems and their applications, \textit{SIAM Review} 23 (1981),
444--466.\bigskip}

This early (1981) review paper brings together and discusses theory and
applications of methods, identified and labelled as row-action methods, for
linear feasibility problems, linearly constrained optimization problems and
some interval convex programming problems. The main feature of row-action
methods is that they are iterative procedures which, without making any
changes to the original matrix $A$, use the rows of $A$, one row at a time.
Such methods are important and have demonstrated effectiveness for problems
with large or huge matrices which do not enjoy any detectable or usable
structural pattern, apart from a high degree of sparseness. Fields of
application where row-action methods are used in various ways include image
reconstruction from projection, operations research and game theory, learning
theory, pattern recognition and transportation theory. A row-action method for
the nonlinear convex feasibility problem is also presented.\textbf{\bigskip}

\cite{censor-adm-84} \textbf{Y. Censor, Iterative methods for the convex
feasibility problem, \textit{Annals of Discrete Mathematics} 20 (1984),
83--91.\bigskip}

Abstract: \textquotedblleft The problem of finding a point in the intersection
of a finite family of closed convex sets in the Euclidean space is considered
here. Several iterative methods for its solution are reviewed and some
connections between them are pointed out.\textquotedblright\textbf{\bigskip}

\cite{coap} \textbf{Y. Censor, W. Chen, P.L. Combettes, R. Davidi and G.T.
Herman, On the effectiveness of projection methods for convex feasibility
problems with linear inequality constraints, \textit{Computational
Optimization and Applications} 51 (2012), 1065--1088.\bigskip}

Besides interesting experimental results, this paper contains (in its
introduction) many pointers to applications wherein projection methods were
used. Section 3 contains a brief glimpse into some recently published results
that show the efficacy of projection methods for some large problems, and
their use in commercial devices.\textbf{\bigskip}

\cite{cen-segal} \textbf{Y. Censor and A. Segal, Iterative projection methods
in biomedical inverse problems, in: Y. Censor, M. Jiang and A.K. Louis
(Editors), \textit{Mathematical Methods in Biomedical Imaging and
Intensity-Modulated Radiation Therapy (IMRT)}, Edizioni della Normale, Pisa,
Italy, 2008, pp. 65--96.\bigskip}

In this paper on projection methods the authors review Bregman projections and
the following algorithmic structures: sequential projection algorithms,
string-averaging algorithmic structures, block-iterative algorithmic schemes
with underrelaxed Bregman projections, component averaging (CAV)).\textbf{ }
Seminorm-induced oblique projections for sparse nonlinear convex feasibility
problems and BICAV (Block-iterative component averaging) are discussed,
followed by a review of subgradient projections and perturbed projections for
the multiple-sets split feasibility problem. Finally, algorithms for the
quasi-convex feasibility problem are presented.\textbf{\bigskip}

\cite{c93} \textbf{P.L. Combettes, The foundations of set-theoretic
estimation,\ \textit{Proceedings of the IEEE} 81 (1993), 182--208.\bigskip}

The paper eloquently presents the rational of posing problems as feasibility
problems rather than optimization problems. Well-connected to the practical
aspects of estimation problems, the discussion starts with a study of various
sets to define solutions and continues with mathematical methods consistent or
inconsistent feasibility problems. A comprehensive section on connections with
other estimation procedures closes this interesting (equipped with 223
references) paper.\textbf{\bigskip}

\cite{c96} \textbf{P.L. Combettes, The convex feasibility problem in image
recovery,\ in: \textit{Advances in Imaging and Electron Physics}, vol.
95,\ (P. Hawkes, Editor), pp. 155--270, Academic Press, New York, NY, USA,
1996.\bigskip}

Image recovery is a broad discipline that encompasses the large body of
inverse problems in which an image is to be inferred from the observation of
data consisting of signals physically or mathematically related to it. Image
restoration and image reconstruction from projections are two main
sub-branches of image recovery. The traditional approach has been to use a
criterion of optimality, which usually leads to a single \textquotedblleft
best\textquotedblright\ solution. An alternative approach is to use
feasibility as an acceptance test, in which compliance with all prior
information and the data defines a set of equally acceptable solutions. This
framework, which underlies the feasibility-seeking approach, is discussed in
this survey. It contains an overview of convex set theoretic image recovery,
of construction of property sets, of solving the convex feasibility problem,
and numerical examples.\textbf{\bigskip}

\cite{comb-pesq} \textbf{P.L. Combettes and J.-C. Pesquet, Proximal splitting
methods in signal processing, in: \textit{Fixed-Point Algorithms for Inverse
Problems in Science and Engineering}, (H.H. Bauschke, R.S. Burachik, P.L.
Combettes, V. Elser, D.R. Luke and H. Wolkowicz, Editors), Springer, New York,
NY, USA, 2011, pp. 185--212.\bigskip}

The proximity operator of a convex function is a natural extension of the
notion of a projection operator onto a convex set. This tool, which plays a
central role in the analysis and the numerical solution of convex optimization
problems, has recently been introduced in the arena of inverse problems and,
especially, in signal processing, where it has become increasingly important.
In this paper, the authors review the basic properties of proximity operators
which are relevant to signal processing and present optimization methods based
on these operators. These proximal splitting methods are shown to capture and
extend several well-known algorithms in a unifying framework, covering several
important projection methods. \textbf{\bigskip}

\cite{dl14} \textbf{F. Deutsch and H. Hundal, Arbitrarily slow convergence of
sequences of linear operators: An updated survey, in: S. Reich and A.
Zaslavski (Editors), Proceedings of the Workshop on Infinite Products of
Operators and Their Applications, Haifa, Israel, 2012, \textit{Contemporary
Mathematics}, accepted for publication.\bigskip}

This is an updated survey on the slowest possible rate of convergence of a
sequence of linear operators that converges pointwise to a linear operator.
Although written in the general context, it contains special sections that
relate the general results to projections: Section 9 is on Application to
Cyclic Projections, and Scetion 10 is on Application to Intermittent
Projections\textbf{. \bigskip}

\cite{youla} \textbf{D.C. Youla, Mathematical theory of image restoration by
the method of convex projections, in: \textit{Image Recovery: Theory and
Applications} (H. Stark, editor), Academic Press, Orlando, FL, USA (1987), pp.
29--78.\bigskip}

Projection methods are frequently called POCS methods (\textquotedblleft
projections onto convex sets\textquotedblright) particularly in the
engineering community. In this 1978 article, Youla suggested that many
problems in image restoration could be formulated in terms of linear subspaces
and orthogonal projections in Hilbert space. He noted that the linear
formulation can result in loss of important information and lead to an
ill-posed restoration problem, and offered orthogonal projections onto closed
convex sets as a way to smooth the restoration problem and reintroduce
information about the image to be restored. Because Youla's objective is to
familiarize the image processing community with the fundamentals of convex
feasibility, asymptotically regular nonexpansive mappings, and fixed point
theory, this paper is a relatively self-contained tutorial introduction to
these topics and to the results known up to 1987.\textbf{\bigskip}

\cite{zarantonnelo71} \textbf{E.H. Zarantonello, Projections on convex sets in
Hilbert space and spectral theory. Part I. Projections on convex sets, in:
\textit{Contributions to Nonlinear Functional Analysis}, E.H. Zarantonello
(Editor), Academic Press, New York, NY, USA, 1971, pp. 239--341.\bigskip}

Everything you wanted probably to know about projections on convex sets in
1971 is probably in this Part I (out of two parts) paper. Its sections are: 1.
Projections, basic properties; 2. Vertices and faces; 3. The range of
$I-P_{K}$; 4.Translation sets, parallel convex sets, differentiability; 5. The
algebra of projections.

\section{The end is open}

In a recent paper \cite{ksp09} it is shown that a variant of ART (i.e., of
Kaczmarz's sequential projections method) can be used for crystal lattice
orientation distribution function estimation from diffraction data. One of the
problems discussed in \cite{ksp09} has 1,372,000,000 unknowns and the number
of equations is potentially infinite. They are randomly generated and a
projection step can be carried out as soon as a new equation is available (an
ideal use of a sequential projection method of the row-action type, see
\cite{yc-siam-81}). The result reported in the paper for that problem is
obtained after 1,000,000,000 such projection steps. As for all methodologies,
projection methods are not necessarily the approach of choice in all
applications. However, in important applications in biomedicine, image
processing, and many other fields, see \cite{coap}, projection methods work
well and have been used successfully for a long time.\bigskip

\textbf{Acknowledgments}. Some of our colleagues to whom we sent an early
version of the paper helped us bring it to its present form. The help, ranging
from encouragements to typos-hunting, included in some cases also suggestions
for addition of specific items that we were unaware of. For this help we are
indebted to Heinz Bauschke, Jon Borwein, Charlie Byrne, John Chinneck, Patrick
Combettes, Frank Deutsch, Tommy Elfving, Masao Fukushima, Gabor Herman, Ming
Jiang, Igor Konnov, Boris Polyak, Constantin Popa, Simeon Reich, Isao Yamada,
and Alexander Zaslavski.


\begin{thebibliography}{99}                                                                                               %


\bibitem {agmon}S. Agmon, The relaxation method for linear inequalities,
\textit{Canadian Journal of Mathematics} \textbf{6} (1954), 382--392.

\bibitem {amemiya-ando}I. Amemiya and T. Ando, Convergence of random products
of contractions in Hilbert space, \textit{Acta Scientiarum Mathematicarum
(Szeged)} \textbf{26} (1965), 239--244.

\bibitem {auslender}A. Auslender, \textit{Optimisation: M\'{e}thodes
Num\'{e}riques}, Masson, Paris, France, 1976\ (in French).

\bibitem {bb96}H.H. Bauschke and J.M. Borwein, On projection algorithms for
solving convex feasibility problems,\ \textit{SIAM Review} \textbf{38} (1996), 367--426.

\bibitem {bbl97}H.H. Bauschke, J.M. Borwein and A.S. Lewis, The method of
cyclic projections for closed convex sets in Hilbert space,
\textit{Contemporary Mathematics} \textbf{204} (1997), 1--38.

\bibitem {bauschke-conf}H.H. Bauschke, R.S. Burachik, P.L. Combettes, V.
Elser, D.R. Luke and H. Wolkowicz (Editors), \textit{Fixed-Point Algorithms
for Inverse Problems in Science and Engineering}, Springer, 2011.

\bibitem {BC11}H.H. Bauschke and P.L. Combettes, \textit{Convex Analysis and
Monotone Operator Theory in Hilbert Spaces}, Springer, New York, NY, USA, 2011.

\bibitem {swiss knives}H.H. Bauschke and V.R. Koch, Projection methods: Swiss
army knives for solving feasibility and best approximation problems with
halfspaces, in: S. Reich and A. Zaslavski (Editors), \textit{Proceedings of
the workshop on Infinite Products of Operators and Their Applications}, Haifa,
Israel, 2012, \textit{Contemporary Mathematics}, accepted for publication.
Available at: http://arxiv.org/abs/1301.4506.

\bibitem {bg03}A. Ben-Israel and T. Greville, \textit{Generalized Inverses:
Theory and Applications}, 2nd Edition, Springer-Verlag, New York, NY, USA, 2003.

\bibitem {benzi05}M. Benzi, Gianfranco Cimmino's contributions to numerical
mathematics, in: \textit{Atti del Seminario di Analisi Matematica,
Dipartimento di Matematica dell'Universita` di Bologna}. Volume Speciale:
Ciclo di Conferenze in Memoria di Gianfranco Cimmino, Marzo-Aprile 2004,
Tecnoprint, Bologna, Italy (2005), pp. 87--109. Available at:
http://www.mathcs.emory.edu/\symbol{126}benzi/Web\_papers/cimmino.pdf.

\bibitem {berinde}V. Berinde, \textit{Iterative Approximation of Fixed
Points}, Lecture Notes in Mathematics 1912, Springer-Verlag, Berlin,
Heidelberg, Germany, 2007.

\bibitem {BO75}E. Blum and W. Oettli, \textit{Mathematische Optimierung:
Grundlagen und Verfahren}, Springer-Verlag, Berlin, Heidelberg, Germany, 1975
(in German).

\bibitem {borwein-zhu-book}J.M. Borwein and Q.J. Zhu, \textit{Techniques of
Variational Analysis}, Springer-Verlag, New York, NY, USA, 2005. Paperback
2010. Available at: http://www.carma.newcastle.edu.au/jon/ToVA/tova-talks.pdf

\bibitem {Bregman}L.M. Bregman, The method of successive projection for
finding a common point of convex sets, \textit{Soviet Mathematics Doklady}
\textbf{6} (1965), 688--692. Translated from the Russian original publication:
\textit{Doklady Akademii Nauk SSSR} \textbf{162} (1965) 487--490.

\bibitem {bregman67}L.M. Bregman, The relaxation method of finding the common
point of convex sets and its application to the solution of problems in convex
programming,\ \textit{USSR Computational Mathematics and Mathematical Physics}
\textbf{7 (}1967), 200--217.

\bibitem {brezinski-book}C. Brezinski, \textit{Projection Methods for Systems
of Equations}, Elsevier Science Publishers, Amsterdam, The Netherlands, 1997.

\bibitem {bruck}R.E. Bruck, On the random product of orthogonal projections in
Hilbert space II, \textit{Contemporary Mathematics} \textbf{513} (2010), 65--98.

\bibitem {pink-book}D. Butnariu, Y. Censor and S. Reich (Editors),
\textit{Inherently Parallel Algorithms in Feasibility and Optimization and
Their Applications}, Elsevier Science Publishers, Amsterdam, The Netherlands, 2001.

\bibitem {byrnebook}C.L. Byrne,\textit{ Applied Iterative Methods, }AK Peters,
Wellsely, MA, USA, 2008.

\bibitem {byrne14book}C.L. Byrne, \textit{Iterative Optimization in Inverse
Problems}, Chapman and Hall/CRC Press, 2014.

\bibitem {CEG12}A. Cegielski, \textit{Iterative Methods for Fixed Point
Problems in Hilbert Spaces}, Lecture Notes in mathematics 2057,
Springer-Verlag, Berlin, Heidelberg, Germany, 2012.

\bibitem {yc-siam-81}Y. Censor, Row-action methods for huge and sparse systems
and their applications, \textit{SIAM Review} \textbf{23} (1981), 444--466.

\bibitem {censor-adm-84}Y. Censor, Iterative methods for the convex
feasibility problem, \textit{Annals of Discrete Mathematics} \textbf{20
}(1984), 83--91.

\bibitem {coap}Y. Censor, W. Chen, P.L. Combettes, R. Davidi and G.T. Herman,
On the effectiveness of projection methods for convex feasibility problems
with linear inequality constraints, \textit{Computational Optimization and
Applications} \textbf{51} (2012), 1065--1088.

\bibitem {cl81}Y. Censor and A. Lent, An iterative row-action method for
interval convex programming,\ \textit{Journal of Optimization Theory and
Applications} \textbf{34 (}1981), 321--353.

\bibitem {cen-segal}Y. Censor and A. Segal, Iterative projection methods in
biomedical inverse problems, in: Y. Censor, M. Jiang and A.K. Louis (Editors),
\textit{Mathematical Methods in Biomedical Imaging and Intensity-Modulated
Radiation Therapy (IMRT)}, Edizioni della Normale, Pisa, Italy, 2008, pp. 65--96.

\bibitem {CZ97}Y. Censor and S.A. Zenios, \textit{Parallel Optimization:
Theory, Algorithms, and Applications}, Oxford University Press, New York, NY,
USA, 1997.

\bibitem {cimmino}G. Cimmino, Calcolo approssimato per soluzioni dei sistemi
di equazioni lineari, \textit{La Ricerca Scientifica XVI}, Series II, Anno IX
\textbf{1} (1938), 326--333 (in Italian).

\bibitem {chinneck-book}J.W. Chinneck, \textit{Feasibility and Infeasibility
in Optimization: Algorithms and Computational Methods,} Springer, New York,
NY, USA, 2007.

\bibitem {c93}P.L. Combettes,\textit{ }The foundations of set-theoretic
estimation,\ \textit{Proceedings of the IEEE} \textbf{81 (}1993), 182--208.

\bibitem {c96}P.L. Combettes, The convex feasibility problem in image
recovery,\ in: \textit{Advances in Imaging and Electron Physics,} vol.
\textbf{95,\ }(P. Hawkes, Editor), pp. 155--270, Academic Press, New York, NY,
USA, 1996.

\bibitem {comb-pesq}P.L. Combettes and J.-C. Pesquet, Proximal splitting
methods in signal processing, in: \textit{Fixed-Point Algorithms for Inverse
Problems in Science and Engineering}, (H.H. Bauschke, R.S. Burachik, P.L.
Combettes, V. Elser, D.R. Luke and H. Wolkowicz, Editors), Springer, New York,
NY, USA, 2011, pp. 185--212.

\bibitem {deutsch01}F. Deutsch, \textit{Best Approximation in Inner Product
Spaces}, Springer-Verlag, New York, NY, USA, 2001.

\bibitem {dl14}F. Deutsch and H. Hundal, Arbitrarily slow convergence of
sequences of linear operators: An updated survey, in: S. Reich and A.
Zaslavski (Editors), \textit{Proceedings of the workshop on Infinite Products
of Operators and Their Applications}, Haifa, Israel, 2012,
\textit{Contemporary Mathematics}, accepted for publication.

\bibitem {eremin65}I.I. Eremin, Generalization of the relaxation method of
Motzkin and Agmon, \textit{Uspekhi Matematicheskikh Nauk} \textbf{20} (1965),
183--187 (in Russian).

\bibitem {eremin}I.I. Eremin, \textit{Theory of Linear Optimization}, VSP --
International Science Publishers, Zeist, The Netherlands, 2002.

\bibitem {ER11}R. Escalante and M. Raydan, \textit{Alternating Projection
Methods}, Society for Industrial and Applied Mathematics (SIAM), Philadelphia,
PA, USA, 2011.

\bibitem {PF03}F. Facchinei and J.-S. Pang, \textit{Finite-Dimensional
Variational Inequalities and Complementarity Problems, Volume I and Volume II,
}Springer-Verlag, New York, NY, USA, 2003.

\bibitem {galantai}A. Gal\'{a}ntai, \textit{Projectors and Projection
Methods}, Kluwer Academic Publishers, Dordrecht, The Netherlands, 2004.

\bibitem {gastinel}N. Gastinel, \textit{Linear Numerical Analysis}, Hermann,
Paris, France, 1970. Translated from the original French text: \textit{Analyse
Num\'{e}rique Lin\'{e}aire, }Hermann, Paris, France, 1966.

\bibitem {gbh71}R. Gordon, R. Bender and G.T. Herman, Algebraic reconstruction
techniques (ART) for three-dimensional electron microscopy and x-ray
photography, \textit{Journal of Theoretical Biology} \textbf{29} (1970), 471--482.

\bibitem {GPR67}L.G. Gubin, B.T. Polyak and E.V. Raik, The method of
projection for finding the common point in convex sets, \textit{USSR
Computational Mathematics and Mathematical Physics} \textbf{7} (1967), 1--24.
Translated from the Russian original publication at: \textit{Zhurnal
Vychislitel'noi Matematiki i Matematicheskoi Fiziki}, \textbf{7} (1967), 1211--1228.

\bibitem {halperin}I. Halperin, The product of projection operators,
\textit{Acta Scientiarum Mathematicarum (Szeged) }\textbf{23} (1962), 96--99.

\bibitem {GTH}G.T. Herman, \textit{Fundamentals of Computerized Tomography:
Image Reconstruction from Projections}, Springer-Verlag, London, UK, 2nd
Edition, 2009.

\bibitem {hounsfield}G.N. Hounsfield, A method and apparatus for examination
of a body by radiation such as X or gamma radiation. UK Patent No. 1283915 (1968/72).

\bibitem {householder64}A.S. Householder, \textit{The Theory of Matrices in
Numerical Analysis}, Dover Publications, Inc., New York, NY, USA, 1975
(originally published by the Blaisdell Publishing Company, a division of Ginn
and Company, New York, NY, USA, 1964).

\bibitem {kacmarcz}S. Kaczmarz, Angen\"{a}herte Aufl\"{o}sung von Systemen
linearer Gleichungen, \textit{Bulletin de l'Acad\'{e}mie Polonaise des
Sciences et Lettres} \textbf{A35} (1937), 355--357 (in German).

\bibitem {kacmarcz-eng93}S. Kaczmarz, Approximate solution of systems of
linear equations, \textit{International Journal of Control} \textbf{57}
(1993), 1269--1271. (Posthumous).

\bibitem {ksp09}I.G. Kazantsev, S. Schmidt and H.F. Poulsen, A discrete
spherical X-ray transform of orientation distribution functions using bounding
cubes, \textit{Inverse Problems} \textbf{25}, 105009 (2009).

\bibitem {konnov}I. Konnov, \textit{Combined Relaxation Methods for
Variational Inequalities}, Lecture Notes in Economics and Mathematical
Systems, Volume 495, Springer-Verlag, Berlin, Heidelberg, Germany, 2001.

\bibitem {kurpel}N.S. Kurpel', \textit{Projection-Iterative Methods for
Solution of Operator Equations}, Translations of Mathematical Monographs
\textbf{46}, American Mathematical Society, Providence, RI, USA, 1976.

\bibitem {maligranda}L. Maligranda, Stefan Kaczmarz (1895-1939),
\textit{Antiquitates Mathematicae} \textbf{1} (2007), 15--61 (in Polish).
Available at: \newline
http://wydawnictwa.ptm.org.pl/index.php/antiquitates-mathematicae. \newline
Electronic reprint in: \textit{Lexicon of Polish Mathematicians} [Leksykon
Matematyk\'{o}w Polskich], 2008, 1--47. Available at:\newline http://leksykon.ptm.mimuw.edu.pl/biogramy/kaczmarz/kaczmarz.php.

\bibitem {mot-schoen54}T.S. Motzkin and I.J. Schoenberg, The relaxation method
for linear inequalities, \textit{Canadian Journal of Mathematics} \textbf{6}
(1954), 393--404.

\bibitem {neumann}J. von Neumann, \textit{Functional Operators, Volume II: The
Geometry of Orthogonal Spaces}, Annals of Mathematics Studies, Volume
\textbf{AM-22}, 1950, Princeton University Press, Princeton, NJ, USA. This is
a reprint of mimeographed lecture notes first distributed in 1933.

\bibitem {parks93}P.C. Parks, Stefan Kaczmarz (1895--1939),
\textit{International Journal of Control} \textbf{57 }(1993), 1263--1267.

\bibitem {p12}C. Popa, \textit{Projection Algorithms - Classical Results and
Developments, Applications to Image Reconstruction}, Lambert Academic
Publishing - AV Akademikerverlag GmbH \& Co. KG, Saarbr\"{u}cken, Germany, 2012.

\bibitem {saad03}Y. Saad, \textit{Iterative Methods for Sparse Linear
Systems}, Second Edition, Society for Industrial and Applied Mathematics
(SIAM), Philadelphia, PA, USA, 2003. Available at: \newline
http://www-users.cs.umn.edu/\symbol{126}saad/IterMethBook\_2ndEd.pdf.

\bibitem {itor}H.D. Scolnik, A.R. De Pierro, N.E. Echebest and M.T.
Guardarucci (Guest Editors), Special Issue of \textit{International
Transactions in Operational Research, }Volume 16, Issue 4, July 2009.

\bibitem {stark}H. Stark and Y. Yang, \textit{Vector Space Projections: A
Numerical Approach to Signal and Image Processing, Neural Nets, and Optics},
John Wiley \& Sons, Inc. New York, NY, USA, 1998.

\bibitem {beckenbach}C.B. Tompkins, Methods of steepest descent, in: E.F.
Beckenbach (Editor), \textit{Modern Mathematics for the Engineer: First
Series}, McGraw-Hill, New York, NY, USA, 1956, pp. 448--479.

\bibitem {youla}D.C. Youla, Mathematical theory of image restoration by the
method of convex projections, in: \textit{Image Recovery: Theory and
Applications }(H. Stark, editor), Academic Press, Orlando FL, USA (1987), pp. 29--78.

\bibitem {zarantonnelo71}E.H. Zarantonello, Projections on convex sets in
Hilbert space and spectral theory. Part I. Projections on convex sets, in:
\textit{Contributions to Nonlinear Functional Analysis}, E.H. Zarantonello
(Editor), Academic Press, New York, NY, USA, 1971, pp. 239--341.
\end{thebibliography}
\end{document}